\documentstyle[11pt,twoside]{article}
\title{\bf A note on a generalisation of a definite integral involving the Bessel function of the first kind}
\author{\sc S.A. Dar$^a$, M. Kamarujjama$^b$ and R.B. Paris$^c$\\
\\
${}^a\!$ {\em Department of Applied Mathematics, Aligarh Muslim University,}\\
{\em Aligarh-202002, India}\\ 
{\em E-Mail: showkatjmi34@gmail.com}\\
${}^b\!$ {\em Department of Applied Mathematics, Aligarh Muslim University }\\
{\em Aligarh-202002, India}\\
{\em E-mail: mdkamarujjama@rediffmail.com}\\ 
${}^c\!$ {\em Division of Computing and Mathematics, Abertay University,}\\
{\em Dundee DD1 1HG, UK}\\
{\em E-Mail: r.paris@abertay.ac.uk}
}
\begin{document}
\newcommand{\bee}{\begin{equation}}
\newcommand{\ee}{\end{equation}}
\def\f#1#2{\mbox{${\textstyle \frac{#1}{#2}}$}}
\def\dfrac#1#2{\displaystyle{\frac{#1}{#2}}}
\newcommand{\fr}{\frac{1}{2}}
\newcommand{\fs}{\f{1}{2}}
\newcommand{\al}{\alpha}
\newcommand{\g}{\Gamma}
\newcommand{\br}{\biggr}
\newcommand{\bl}{\biggl}
\date{}
\maketitle
\pagestyle{myheadings}
\markboth{\hfill \it S.A. Dar, M. Kamarujjama and R.B. Paris  \hfill}
{\hfill \it A note on a generalisation of an integral \hfill}
\begin{abstract}
We consider a generalisation of a definite integral involving the Bessel function of the first kind. It is shown that this integral can be expressed in terms of the Fox-Wright function ${}_p\Psi_q(z)$ of one variable. Some consequences 
of this representation are explored by suitable choice of parameters. In addition, two closed-form evaluations of infinite series of the Fox-Wright function are deduced.
\vspace{0.4cm}

\noindent {\bf 2020 AMS Classification:} 33C20, 33B10, 42A38, 33C60, 44A20
\vspace{0.3cm}

\end{abstract}

\vspace{0.3cm}

\noindent $\,$\hrulefill $\,$

\vspace{0.2cm}

\begin{center}
{\bf 1. \  Introduction and Preliminaries}
\end{center}
\setcounter{section}{1}
\setcounter{equation}{0}
\renewcommand{\theequation}{\arabic{section}.\arabic{equation}}
The Fox-Wright function ${}_p\Psi_q(z)$ of one variable \cite{K1, K2} is given by
\[{}_p\Psi_q\bl[\begin{array}{c} (\al_1, A_1), \ldots ,(\al_p,A_p)\\ (\beta_1, B_1), \ldots ,(\beta_q,B_q)\end{array}\bl|z\br]={}_p\Psi_q\bl[\begin{array}{c} (\al_j, A_j)_{1,p}\\ (\beta_j, B_j)_{1,q} \end{array}\bl|z\br]\]
\bee\label{e11}
=\sum_{k=0}^\infty \frac{\g(\al_1+kA_1) \ldots \g(\al_p+kA_p)}{\g(\beta_1+kB_1)\ldots \g(\beta_q+kB_q)}\,\frac{z^k}{k!},
\ee
where $z\in{\bf C}$, $\al_j$, $\beta_j\in{\bf C}$ and the coefficients $A_j\geq0$, $B_j\geq0$, it being supposed throughout that the $\al_j$, $\beta_j$ and the $A_j$, $B_j$ are such that the gamma functions are well defined. With
\[\Delta=\sum_{j=1}^qB_j-\sum_{j=1}^pA_j,\quad \delta=\prod_{j=1}^p A_j^{-A_j} \prod_{j=1}^qB_j^{B_j},\quad \mu^*=\sum_{j=1}^q\beta_j-\sum_{j=1}^p\al_j+\fs(p-q),\]
the series in (\ref{e11}) converges for $|z|<\infty$ when $\Delta>-1$, for $|z|<\delta$ when $\Delta=-1$  and for $|z|=\delta$ if, in addition, $\Re (\mu)>\fs$.

Particular cases of (\ref{e11}) that we shall employ are the Wright function \cite[p.~438(1.2)]{K2} is defined by 
\bee\label{e12}
\Phi(\al,\beta,z)={}_0\Psi_1\bl[\begin{array}{c}-\!\!-\\(\beta,\al)\end{array}\bl|z\br]=\sum_{k=0}^\infty\frac{z^k}{k! \g(\beta+\al k)},
\ee
where $z, \beta\in{\bf C}$ and $\al>0$ and the Wright generalised-Bessel function \cite[p.~438(1.3)]{K2} defined by
\bee\label{e13}
J_\nu^\mu(-z)=\Phi(\mu,\nu+1,-z)={}_0\Psi_1\bl[\begin{array}{c}-\!-\\(\nu+1,\mu)\end{array}\bl|-z\br]=\sum_{k=0}^\infty\frac{(-z)^k}{k! \g(\nu+1+\mu k)},
\ee
where $z, \nu\in{\bf C}$ and $\mu>0$. The Mittag-Leffler function \cite[p.~450(6.1)]{K2} is given by
\bee\label{e14}
E_{\al,\beta}(z)={}_1\Psi_1\bl[\begin{array}{c} (1,1)\\(\beta,\al)\end{array}\bl|z\br]=\sum_{k=0}^\infty \frac{z^k}{\g(\beta+\al k)},
\ee
where $z, \beta\in{\bf C}$ and $\al>0$.
The Bessel function of the first kind is defined by (see \cite[p.~217]{DLMF})
\[J_\nu(z)=(\fs z)^\nu \sum_{k=0}^\infty\frac{ (-1)^k (z/2)^{2k}}{k! \g(\nu+1+k)}.\]
Of special interest in this paper are the elementary functions corresponding to $\nu=\pm\fs$, namely
\bee\label{e15}
J_{1/2}(z)=\sqrt{\frac{2}{\pi z}}\,\sin z,\qquad J_{-1/2}(z)=\sqrt{\frac{2}{\pi z}}\,\cos z.
\ee

The infinite Fourier sine and cosine transforms of $g(x)$ over the interval $[0,\infty)$ are defined by \cite{Erd, Oberh1}
\bee\label{e16}
F_{S,C}(g(x);b)=\int_0^\infty g(x)\begin{array}{c}\sin\\\cos\end{array}(bx)\,dx=G_{S,C}(b)\qquad (b>0).
\ee
If $\Re(z)>0$ and $\Re (b)>0$, the Mellin transforms of $\cos (bx)/(e^{ax}-1)$ and $\sin (bx)/(e^{ax}-1)$ are given by \cite[p.~15, 1.4(7)]{Erd}
\[\int_0^\infty \frac{x^{\mu-1} \cos (bx)}{e^{ax}-1}\,dx=\frac{\g(\mu)}{2a^\mu}\bl\{\zeta(\mu,1+ib/a)+\zeta(\mu,1-ib/a)\br\}\]
and
\[\int_0^\infty \frac{x^{\mu-1} \sin (bx)}{e^{ax}-1}\,dx=\frac{i\g(\mu)}{2a^\mu}\bl\{\zeta(\mu,1+ib/a)-\zeta(\mu,1-ib/a)\br\},\]
where $\zeta(a,z)=\sum_{k\geq0}(k+z)^{-a}$ is the Hurwitz zeta function. In the special case $\mu=2$ and $a=2\pi$, $b=\pi n$, we have
\bee\label{e17a}
\int_0^\infty \frac{x \cos (\pi nx)}{e^{2\pi x}-1}dx=\frac{1}{2\pi^2n^2}+\frac{1}{4(1-\cosh \pi n)}
\ee
and \cite[(5.15.1)]{DLMF}
\bee\label{e17b}
\int_0^\infty \frac{x \sin (\pi nx)}{e^{2\pi x}-1}dx=\frac{i}{8\pi^2}\{\psi'(1+\fs in)-\psi'(1-\fs in)\},
\ee
where $\psi(z)$ is the logarithmic derivative of $\g(z)$ and the prime denotes differentiation.

A natural generalisation of the Gauss hypergeometric function ${}_2F_1(z)$ is the generalised hypergeometric function ${}_pF_q(z)$, with $p$ numerator parameters $\al_1, \ldots , \al_p$ and $q$ denominator parameters $\beta_1, \ldots  ,\beta_q$, defined by
\bee\label{e18}
{}_pF_q\bl(\begin{array}{c} \al_1, \ldots , \al_p\\\beta_1, \ldots  , \beta_q\end{array};z\br)=
{}_pF_q\bl(\begin{array}{c} (\al_j)_{1,p}\\(\beta_j)_{1,q}\end{array};z\br)=
\sum_{k=0}^\infty \frac{(\al_1)_k \ldots (\al_p)_k}{(\beta_1)_k \ldots (\beta_q)_k}\,\frac{z^k}{k!},
\ee
where $\al_j\in{\bf C}$ ($j=1, \ldots ,p$) and $\beta_j\in{\bf C}\backslash {\bf Z}_0^-$ ($j=1, \ldots ,q$), ${\bf Z}_0^-=\{0,-1, -2, \ldots \}$ and $(a)_k=\g(a+k)/\g(a)$ is the Pochhammer symbol. The series in (\ref{e18}) is convergent for $|z|<\infty$ if $p\leq q$ and for $|z|<1$ if $p=q+1$. If we define the parametric excess by $\omega=\sum_{j=1}^q\beta_j-\sum_{j=1}^p\al_j$, then it is known that the ${}_pF_q(z)$ series with $p=q+1$ is (i) absolutely convergent for $|z|=1$ if $\Re (\omega)>0$, (ii) is conditionally converegnt for $|z|=1$ if $-1<\Re (\omega)\leq 0$ and (iii) is divergent for $|z|=1$ if $\Re (\omega)\leq -1$.

The central aim in this paper is to give a generalisation and find some consequences of a definite integral involving the Bessel function of the first kind which we express in terms of the Fox-Wright $\Psi$ function. The work is motivated by the papers by one of the present authors in \cite{Dar1, Dar2, Dar3}. In order to generalise the definite integrals  introduced by Ramanujan $\phi_{S,C}(m,n)$ defined by
\bee\label{e19}
\phi_{S,C}(m,n)=\int_0^\infty \frac{x^m}{e^{2\pi x}-1}\,\begin{array}{c}\sin\\\cos\end{array}(\pi nx)\,dx,
\ee
we introduce the following integrals:
\bee\label{e110}
F_1(\mu,\xi,a,\nu,y)=\int_0^\infty x^\mu e^{-ax^\xi} \sqrt{xy}\,J_\nu(xy)\,dx,
\ee
\bee\label{e111}
F_2(\mu,\xi,b,c,\nu,y)=\sum_{k=0}^\infty\frac{\Theta(k)}{k!} \int_0^\infty x^\mu e^{-(b+ck)x^\xi} \sqrt{xy}\,J_\nu(xy)\,dx,
\ee
\bee\label{e112}
F_3(\mu,\xi,b,c,\nu,y)=\int_0^\infty x^\mu e^{-bx^\xi} {}_r\Psi_s\bl[\begin{array}{c}(\al_j,A_j)_{1,r}\\(\beta_j,B_j)_{1,s}\end{array}\bl|e^{-cx^\xi}\br] \sqrt{xy}\,J_\nu(xy)\,dx
\ee
and
\bee\label{e113}
F_4(\mu,\xi,b,c,\nu,y)=\int_0^\infty x^\mu e^{-bx^\xi} {}_rF_s\bl(\begin{array}{c}(\al_j)_{1,r}\\(\beta_j)_{1,s}\end{array};e^{-cx^\xi}\br)\sqrt{xy}\,J_\nu(xy)\,dx.
\ee 
Here $\xi>0$ and $\{\Theta(k)\}_{k=0}^\infty$ is a bounded sequence of real or complex quantities. We show how the main general theorem given in Section 3 is applicable for obtaining new and interesting results by suitable adjustment of the parameters and variables.
\vspace{0.6cm}

\begin{center}
{\bf 2. \  Evaluation of the definite integral $F_1(\mu,\xi,a,\nu,y)$}
\end{center}
\setcounter{section}{2}
\setcounter{equation}{0}
\renewcommand{\theequation}{\arabic{section}.\arabic{equation}}
In this section, we evaluate the integral in (\ref{e110}) involving the Bessel function of the first kind in terms of the Fox-Wright function. We suppose throughout that the parameter $\xi>0$. 
\newtheorem{theorem}{Theorem}
\begin{theorem}$\!\!\!.$\ \ 
With $\sigma:=(2\mu+2\nu+3)/(2\xi)$, we have
\bee\label{e21}
F_1(\mu,\xi,a,\nu,y)=\int_0^\infty x^\mu e^{-ax^\xi} \sqrt{xy}\,J_\nu(xy)\,dx=
\frac{y^{\nu+1/2}}{2^\nu\xi a^\sigma}\,{}_1\Psi_1\bl[\begin{array}{c} (\sigma, 2/\xi)\\(\nu+1,1)\end{array}\bl|-\frac{y^2}{4a^{2/\xi}}\br],
\ee
where $a>0$, $y>0$ and $\Re (\mu+\nu)>-\fs$.
\end{theorem}

\noindent{\it Proof:} \ \ Expanding the Bessel function in its series form, followed by reversal of the order of summation and integration, we find
\[F_1(\mu,\xi,a,\nu,y)=\frac{y^{\nu+1/2}}{2^\nu \xi} \sum_{k=0}^\infty \frac{(-1)^k (y^2/4)^k}{k! \g(\nu+1+k)} \int_0^\infty t^{\sigma+2k/\xi-1} e^{-at}dt\]
\[=\frac{y^{\nu+1/2}}{2^\nu\xi a^\sigma} \sum_{k=0}^\infty \frac{\g(\sigma+2k/\xi)}{k! \g(\nu+1+k)} \bl(\frac{-y^2}{4a^{2/\xi}}\br)^k\]
upon evaluation of the integral as a gamma function, 
where $\sigma=(2\mu+2\nu+3)/(2\xi)$. If we now employ the definition of the Fox-Wright function in (\ref{e11}) in the above series, we obtain the right-hand side of (\ref{e21}).

If we set $\xi=1$ and $\nu=\pm\fs$ in (\ref{e21}), we obtain the Fourier sine and cosine transforms of $x^{\eta-1} e^{-ax}$ in the form 
\[\int_0^\infty x^{\eta-1} e^{-ax} \cos (xy)\,dx=\frac{\g(\eta)}{a^\eta} \sum_{k=0}^\infty \frac{(\eta)_{2k}}{(\fs)_k k!}\,\bl(-\frac{y^2}{4a^2}\br)^k\hspace{2cm}\]
\[\hspace{3.2cm}=\frac{\g(\eta)}{a^\eta}\ {}_2F_1\bl(\begin{array}{c}\fs\eta,\fs\eta+\fs\\ \fs\end{array};-\frac{y^2}{a^2}\br)\]
and
\[\int_0^\infty x^{\eta-1} e^{-ax} \sin (xy)\,dx=\frac{y\g(\eta+1)}{a^{\eta+1}} \sum_{k=0}^\infty \frac{(\eta+1)_{2k}}{(\f{3}{2})_k k!}\,\bl(-\frac{y^2}{4a^2}\br)^k\hspace{3cm}\]
\[\hspace{2.5cm}=\frac{y\g(\eta+1)}{a^{\eta+1}}\ {}_2F_1\bl(\begin{array}{c}\fs\eta+\fs,\fs\eta+1\\\f{3}{2}\end{array};-\frac{y^2}{a^2}\br).\]
Use of the standard evaluations of the hypergeometric function given in \cite[(15.4.8), (15.4.10)]{DLMF}
then yields the results stated in \cite{Oberh1}
\[\int_0^\infty x^{\eta-1} e^{-ax} \begin{array}{c}\sin\\\cos\end{array} (xy)\,dx=\g(\eta)\, (a^2+y^2)^{-\eta/2}\begin{array}{c}\sin\\\cos\end{array}(\eta \arctan (y/a)),\]
where $\Re (a)>0$, $y>0$ and $\Re (\eta)>-1$ and $\Re (\eta)>0$ for the sine and cosine integral, respectively.
\vspace{0.6cm}

\begin{center}
{\bf 3. \  Evaluation of the definite integrals $F_j(\mu,\xi,a,\nu,y)$, $j=2,3,4$}
\end{center}
\setcounter{section}{3}
\setcounter{equation}{0}
\renewcommand{\theequation}{\arabic{section}.\arabic{equation}}
Here we let $\sigma=(2\mu+2\nu+3)/(2\xi)$ throughout this section. 
\begin{theorem}$\!\!\!.$\ \ Let $\{\Theta(k)\}_{k=0}^\infty$ be a bounded sequence of arbitrary real or complex numbers. Then when $\xi>0$ and $\Re (b)>0$ we have
\[F_2(\mu,\xi,b,c,\nu,y)=\sum_{k=0}^\infty\frac{\Theta(k)}{k!} \int_0^\infty x^\mu e^{-(b+ck)x^\xi} \sqrt{xy}\,J_\nu(xy)\,dx\hspace{3cm}\]
\bee\label{e31}
\hspace{2cm}=\frac{y^{\nu+1/2}}{2^\nu b^\sigma\xi} \sum_{k=0}^\infty \frac{\Theta(k)}{k!} \bl(\frac{(b/c)_k}{(1+b/c)_k}\br)^\sigma
\,{}_1\Psi_1\bl[\begin{array}{c}(\sigma,2/\xi)\\(\nu+1,1)\end{array}\bl|-\frac{y^2}{4(b+ck)^{2/\xi}}\br].\ \ 
\ee
\end{theorem}

\noindent{\it Proof:}\ \ The proof follows the same procedure as in Theorem 1 by expressing the Bessel function in series form and integrating term by term. We find
\[F_2(\mu,\xi,b,c,\nu,y)=\frac{y^{\nu+1/2}}{2^\nu \xi} \sum_{k=0}^\infty \frac{\Theta(k)}{k! (b+ck)^\sigma}
\sum_{\ell=0}^\infty \frac{\g(\sigma+2\ell/\xi)}{\ell! \g(\nu+1+\ell)}\,\bl(\frac{-y^2}{4(b+ck)^{2/\xi}}\br)^\ell.\]
Employing the definition of the Fox-Wright function function in (\ref{e11}) then yields the stated result in (\ref{e31}).
\bigskip

\noindent{\bf Corollary 1.}\ \ For $\xi>0$, $\Re (b)>0$ and $\Re (c)>0$, we have
\[F_3(\mu,\xi,b,c,\nu,y)=\int_0^\infty x^\mu e^{-bx^\xi} {}_r\Psi_s\bl[\begin{array}{c}(\al_j,A_j)_{1,r}\\(\beta_j,B_j)_{1,s}\end{array}\bl|e^{-cx^\xi}\br] \sqrt{xy}\,J_\nu(xy)\,dx\]
\bee\label{e32}
=\frac{y^{\nu+1/2}}{2^\nu b^\sigma\xi} \sum_{k=0}^\infty \frac{\prod_{j=1}^r\g(\al_j+kA_j)}{k! \prod_{j=1}^s\g(\beta_j+kB_j)} \bl(\frac{(b/c)_k}{(1+b/c)_k}\br)^\sigma
\,{}_1\Psi_1\bl[\begin{array}{c}(\sigma,2/\xi)\\(\nu+1,1)\end{array}\bl|-\frac{y^2}{4(b+ck)^{2/\xi}}\br],
\ee
where the parameters $\al_j, \beta_j\in{\bf C}$ and $A_j\geq0$, $B_j\geq0$. This result follows immediately from (\ref{e31}) by substituting
\[\Theta(k)=\frac{\g(\al_1+kA_1)\ldots \g(\al_r+kA_r)}{\g(\beta_1+kB_1)\ldots \g(\beta_s+kB_s)}\qquad (k=0,1,2, \ldots).\]
\vspace{0.1cm}

\noindent{\bf Corollary 2.}\ \ For $\xi>0$, $\Re (b)>0$ and $\Re (c)>0$, we have
\[F_4(\mu,\xi,b,c,\nu,y)=\int_0^\infty x^\mu e^{-bx^\xi} {}_rF_s\bl(\begin{array}{c}(\al_j)_{1,r}\\(\beta_j)_{1,s}\end{array};e^{-cx^\xi}\br)\sqrt{xy}\,J_\nu(xy)\,dx\]
\bee\label{e33}
=\frac{y^{\nu+1/2}}{2^\nu b^\sigma\xi} \sum_{k=0}^\infty \frac{\prod_{j=1}^r(\al_j)_k}{k! \prod_{j=1}^s(\beta_j)_k} \bl(\frac{(b/c)_k}{(1+b/c)_k}\br)^\sigma
\,{}_1\Psi_1\bl[\begin{array}{c}(\sigma,2/\xi)\\(\nu+1,1)\end{array}\bl|-\frac{y^2}{4(b+ck)^{2/\xi}}\br],
\ee
where the parameters $\al_j\in{\bf C}$, $\beta_j\in{\bf C}\backslash{\bf Z}_0^-$ and $r\leq s+1$. This result follows by setting $A_1=\ldots =A_r=1$ and $B_1=\ldots =B_s=1$ in (\ref{e32}).
\vspace{0.2cm}

\noindent{\bf 3.1\ \ Special cases of the integral (\ref{e32})}
\vspace{0.1cm}

\noindent
Special cases of the integral in (\ref{e32}) are given by the following.
When $r=0$, $s=1$ in (\ref{e32}), we obtain
\[\int_0^\infty x^\mu e^{-bx^\xi} \Phi(B_1, \beta_1, e^{-cx^\xi}) \sqrt{xy}\,J_\nu(xy)\,dx\]
\bee\label{e34}
=\frac{y^{\nu+1/2}}{2^\nu b^\sigma\xi} \sum_{k=0}^\infty \frac{1}{k! \g(\beta_1+kB_1)} \bl(\frac{(b/c)_k}{(1+b/c)_k}\br)^\sigma\,{}_1\Psi_1\bl[\begin{array}{c}(\sigma, 2/\xi)\\(\nu+1,1)\end{array}\bl| -\frac{y^2}{4(b+ck)^{2/\xi}}\br],
\ee
where $\beta_1\in{\bf C}$, $B_1>0$ and $\Phi(B_1, \beta_1, z)$ is the Wright function defined in (\ref{e12}).

When $r=0$, $s=1$, $B_1=\mu$, $\beta_1=\gamma+1$ and $e^{-cx^\xi}$ is replaced by $-e^{-cx^\xi}$ in (\ref{e32}), we obtain
\[\int_0^\infty x^\mu e^{-bx^\xi} J_\gamma^\mu(-e^{-cx^\xi})\,\sqrt{xy}\,J_\nu(xy)\,dx\]
\bee\label{e35}
=\frac{y^{\nu+1/2}}{2^\nu b^\sigma\xi} \sum_{k=0}^\infty \frac{1}{k! \g(\gamma+1+\mu k)} \bl(\frac{(b/c)_k}{(1+b/c)_k}\br)^\sigma\,{}_1\Psi_1\bl[\begin{array}{c}(\sigma, 2/\xi)\\(\nu+1,1)\end{array}\bl| -\frac{y^2}{4(b+ck)^{2/\xi}}\br],
\ee
where $\gamma\in{\bf C}$, $\mu>0$ and $J_\gamma^\mu(z)$ is the Wright generalised Bessel function defined in (\ref{e13}).

When $r=s=1$ and $\al_1=A_1=1$ in (\ref{e32}), we obtain
\[\int_0^\infty x^\mu e^{-bx^\xi} E_{B_1, \beta_1}(e^{-cx^\xi})\,\sqrt{xy}\,J_\nu(xy)\,dx\]
\bee\label{e36}
=\frac{y^{\nu+1/2}}{2^\nu b^\sigma\xi} \sum_{k=0}^\infty \frac{1}{k! \g(\beta_1+kB_1)} \bl(\frac{(b/c)_k}{(1+b/c)_k}\br)^\sigma\,{}_1\Psi_1\bl[\begin{array}{c}(\sigma, 2/\xi)\\(\nu+1,1)\end{array}\bl| -\frac{y^2}{4(b+ck)^{2/\xi}}\br],
\ee
where $\beta_1\in{\bf C}$, $B_1>0$ and $E_{B_1,\beta_1}(z)$ is the Mittag-Leffler function defined in (\ref{e14}).

\vspace{0.6cm}

\begin{center}
{\bf 4. \  Expressions for Fourier cosine and sine transforms}
\end{center}
\setcounter{section}{4}
\setcounter{equation}{0}
\renewcommand{\theequation}{\arabic{section}.\arabic{equation}}
The expressions for the Fourier cosine and sine transforms of $x^{\eta-1} e^{-(b+ck)x}$ in terms of infinite series of the Fox-Wright function ${}_1\Psi_1(\cdot)$ will hold true for $y>0$, $\Re (\eta)>0$, $\Re (b)>0$ and $\Re (c)>0$. If we let $\mu=\eta-1$, $\xi=1$ and $\nu=\pm\fs$ in the main theorem in (\ref{e31}), we obtain after some simplification the following expressions:

\[F_2(\eta-1,1,b,c,-\fs,y)\equiv {\bf F}_C^{(1)}(\eta,b,c,y)=\sqrt{\frac{2}{\pi}}\sum_{k=0}^\infty \frac{\Theta(k)}{k!} \int_0^\infty x^{\eta-1}
e^{-(b+ck)x} \cos (xy)\,dx\]
\bee\label{e41}
=\sqrt{2} \sum_{k=0}^\infty \frac{\Theta(k)}{k!}\,\frac{1}{(b+ck)^\eta} \,{}_1\Psi_1\bl[\begin{array}{c} (\eta,2)\\(\fs,1)\end{array}\bl|-\frac{y^2}{4(b+ck)^2}\br],
\ee

\[F_2(\eta-1,1,b,c,\fs,y)\equiv {\bf F}_S^{(1)}(\eta,b,c,y)=\sqrt{\frac{2}{\pi}}\sum_{k=0}^\infty \frac{\Theta(k)}{k!} \int_0^\infty x^{\eta-1}
e^{-(b+ck)x} \sin (xy)\,dx\]
\bee\label{e42}
=\frac{y}{\sqrt{\pi}} \sum_{k=0}^\infty \frac{\Theta(k)}{k!}\,\frac{1}{(b+ck)^{\eta+1}} \,{}_1\Psi_1\bl[\begin{array}{c} (\eta+1,2)\\(\f{3}{2},1)\end{array}\bl|-\frac{y^2}{4(b+ck)^2}\br],
\ee

\[F_3(\eta-1,1,b,c,-\fs,y)\equiv {\bf F}_C^{(2)}(\eta,b,c,y)\hspace{6cm}\]
\[=\sqrt{\frac{2}{\pi}}\int_0^\infty x^{\eta-1} e^{-bx} \,{}_r\Psi_s\bl[\begin{array}{c}(\al_j,A_j)_{1,r}\\(\beta_j,B_j)_{1,s}\end{array}\bl| e^{-cx}\br]\,\cos (xy)\,dx\hspace{2cm}\]
\bee\label{e43}
=\sqrt{2} \sum_{k=0}^\infty \frac{\prod_{j=1}^r\g(\al_j+kA_j)}{k! \prod_{j=1}^s\g(\beta_j+kB_j)}\,\frac{1}{(b+ck)^\eta} \,{}_1\Psi_1\bl[\begin{array}{c} (\eta,2)\\(\fs,1)\end{array}\bl|-\frac{y^2}{4(b+ck)^2}\br],
\ee

\[F_3(\eta-1,1,b,c,\fs,y)\equiv {\bf F}_S^{(2)}(\eta,b,c,y)\hspace{6cm}\]
\[=\sqrt{\frac{2}{\pi}}\int_0^\infty x^{\eta-1} e^{-bx} \,{}_r\Psi_s\bl[\begin{array}{c}(\al_j,A_j)_{1,r}\\(\beta_j,B_j)_{1,s}\end{array}\bl| e^{-cx}\br]\,\sin (xy)\,dx\hspace{2.6cm}\]
\bee\label{e44}
\hspace{1.4cm}=\frac{y}{\sqrt{\pi}} \sum_{k=0}^\infty \frac{\prod_{j=1}^r\g(\al_j+kA_j)}{k! \prod_{j=1}^s\g(\beta_j+kB_j)}\,\frac{1}{(b+ck)^{\eta+1}} \,{}_1\Psi_1\bl[\begin{array}{c} (\eta+1,2)\\(\f{3}{2},1)\end{array}\bl|-\frac{y^2}{4(b+ck)^2}\br],
\ee

\[F_4(\eta-1,1,b,c,-\fs,y)\equiv {\bf F}_C^{(3)}(\eta,b,c,y)=\sqrt{\frac{2}{\pi}}\int_0^\infty x^{\eta-1} e^{-bx} \,{}_rF_s\bl[\begin{array}{c}(\al_j)_{1,r}\\(\beta_j)_{1,s}\end{array}\bl| e^{-cx}\br]\,\cos (xy)\,dx\]
\bee\label{e45}
=\sqrt{2} \sum_{k=0}^\infty \frac{\prod_{j=1}^r\g(\al_j)_k}{k! \prod_{j=1}^s\g(\beta_j)_k}\,\frac{1}{(b+ck)^\eta} \,{}_1\Psi_1\bl[\begin{array}{c} (\eta,2)\\(\fs,1)\end{array}\bl|-\frac{y^2}{4(b+ck)^2}\br],
\ee

\[F_4(\eta-1,1,b,c,\fs,y)\equiv {\bf F}_S^{(3)}(\eta,b,c,y)=\sqrt{\frac{2}{\pi}}\int_0^\infty x^{\eta-1} e^{-bx} \,{}_rF_s\bl[\begin{array}{c}(\al_j)_{1,r}\\(\beta_j)_{1,s}\end{array}\bl| e^{-cx}\br]\,\sin (xy)\,dx\]
\bee\label{e46}
=\frac{y}{\sqrt{\pi}} \sum_{k=0}^\infty \frac{\prod_{j=1}^r\g(\al_j)_k}{k! \prod_{j=1}^s\g(\beta_j)_k}\,\frac{1}{(b+ck)^{\eta+1}} \,{}_1\Psi_1\bl[\begin{array}{c} (\eta+1,2)\\(\f{3}{2},1)\end{array}\bl|-\frac{y^2}{4(b+ck)^2}\br].
\ee
\vspace{0.6cm}

\begin{center}
{\bf 5. \  Closed-form infinite summation formulas}
\end{center}
\setcounter{section}{5}
\setcounter{equation}{0}
\renewcommand{\theequation}{\arabic{section}.\arabic{equation}}
If we choose $\Theta(k)=k!$, $\eta=m+1$, $y=\pi n$ and $b=c=2\pi$ in (\ref{e41})--(\ref{e46}), we find the following representations involving infinite sums of the Fox-Wright function:
\[\int_0^\infty \frac{x^m \cos (\pi nx)}{e^{2\pi x}-1}\,dx=\sum_{k=0}^\infty \int_0^\infty x^m e^{-2\pi x(1+k)} \cos (\pi nx)\,dx\hspace{5cm}\]
\bee\label{e51}
=\frac{\sqrt{\pi}}{(2\pi)^{m+1}} \sum_{k=0}^\infty \frac{1}{(k+1)^{m+1}}\,{}_1\Psi_1\bl[\begin{array}{c}(m+1,2)\\(\fs,1)\end{array}\bl|-\frac{n^2}{16(1+k)^2}\br],
\ee

\[\int_0^\infty \frac{x^m \sin (\pi nx)}{e^{2\pi x}-1}\,dx=\sum_{k=0}^\infty \int_0^\infty x^m e^{-2\pi x(1+k)} \sin (\pi nx)\,dx\hspace{5cm}\]
\bee\label{e52}
=\frac{\sqrt{\pi}\, n}{4(2\pi)^{m+1}} \sum_{k=0}^\infty \frac{1}{(k+1)^{m+2}}\,{}_1\Psi_1\bl[\begin{array}{c}(m+2,2)\\(\f{3}{2},1)\end{array}\bl|-\frac{n^2}{16(1+k)^2}\br],
\ee

\[\int_0^\infty x^m e^{-2\pi x} {}_r\Psi_s\bl[\begin{array}{c} (\al_j,A_j)_{1,r}\\(\beta_j,B_j)_{1,s}\end{array}\bl| e^{-2\pi x}\br]\,\cos (\pi nx)\,dx\hspace{5cm}\]
\bee\label{e53}
= \frac{\sqrt{\pi}}{(2\pi)^{m+1}}\sum_{k=0}^\infty\frac{\prod_{j=1}^r\g(\al_j+kA_j)}{\prod_{j=1}^s\g(\beta_j+kB_j)}\, \frac{1}{(k+1)^{m+1}}\,{}_1\Psi_1\bl[\begin{array}{c}(m+1,2)\\(\fs,1)\end{array}\bl|-\frac{n^2}{16(1+k)^2}\br],
\ee

\[\int_0^\infty x^m e^{-2\pi x} {}_r\Psi_s\bl[\begin{array}{c} (\al_j,A_j)_{1,r}\\(\beta_j,B_j)_{1,s}\end{array}\bl| e^{-2\pi x}\br]\,\sin (\pi nx)\,dx\hspace{5cm}\]
\bee\label{e54}
=\frac{\sqrt{\pi}\, n}{4(2\pi)^{m+1}}  \sum_{k=0}^\infty\frac{\prod_{j=1}^r\g(\al_j+kA_j)}{\prod_{j=1}^s\g(\beta_j+kB_j)}\, \frac{1}{(k+1)^{m+2}}\,{}_1\Psi_1\bl[\begin{array}{c}(m+2,2)\\(\f{3}{2},1)\end{array}\bl|-\frac{n^2}{16(1+k)^2}\br],
\ee

\[\int_0^\infty x^m e^{-2\pi x} {}_rF_s\bl[\begin{array}{c} (\al_j)_{1,r}\\(\beta_j)_{1,s}\end{array}\bl| e^{-2\pi x}\br]\,\cos (\pi nx)\,dx\hspace{5cm}\]
\bee\label{e55}
=\frac{\sqrt{\pi}}{(2\pi)^{m+1}} \sum_{k=0}^\infty\frac{\prod_{j=1}^r\g(\al_j)_k}{\prod_{j=1}^s\g(\beta_j)_k}\, \frac{1}{(k+1)^{m+1}}\,{}_1\Psi_1\bl[\begin{array}{c}(m+1,2)\\(\fs,1)\end{array}\bl|-\frac{n^2}{16(1+k)^2}\br],
\ee

\[\int_0^\infty x^m e^{-2\pi x} {}_rF_s\bl[\begin{array}{c} (\al_j,A_j)_{1,r}\\(\beta_j,B_j)_{1,s}\end{array}\bl| e^{-2\pi x}\br]\,\sin (\pi nx)\,dx\hspace{5cm}\]
\bee\label{e56}
=\frac{\sqrt{\pi}\, n}{4(2\pi)^{m+1}} \sum_{k=0}^\infty\frac{\prod_{j=1}^r\g(\al_j)_k}{\prod_{j=1}^s\g(\beta_j)_k}\, \frac{1}{(k+1)^{m+2}}\,{}_1\Psi_1\bl[\begin{array}{c}(m+2,2)\\(\f{3}{2},1)\end{array}\bl|-\frac{n^2}{16(1+k)^2}\br].
\ee

Finally, from (\ref{e17a}), (\ref{e17b}), (\ref{e51}) and (\ref{e52}), we obtain the summation formulas
\bee\label{e57}
\sum_{k=0}^\infty \frac{1}{(1+k)^2}\,{}_1\Psi_1\bl[\begin{array}{c} (2,2)\\(\fs,1)\end{array}\bl|-\frac{n^2}{16(1+k)^2}\br]=\sqrt{\pi} \bl(\frac{2}{\pi n^2}+\frac{\pi}{1-\cosh \pi n}\br),
\ee
and
\bee\label{e58}
\sum_{k=0}^\infty \frac{1}{(1+k)^3}\,{}_1\Psi_1\bl[\begin{array}{c} (3,2)\\(\f{3}{2},1)\end{array}\bl|-\frac{n^2}{16(1+k)^2}\br]= \frac{2i}{\sqrt{\pi}\,n} \bl\{\psi'(1+\fs in)-\psi'(1-\fs in)\br\},
\ee
where $n$ is free to be chosen. These sums have been verified numerically using {\it Mathematica}.
\vspace{0.6cm}

\begin{center}
{\bf 6. \  Conclusion}
\end{center}
\setcounter{section}{6}
\setcounter{equation}{0}
\renewcommand{\theequation}{\arabic{section}.\arabic{equation}}
In this note we have shown that a certain integral involving the Bessel function of the first kind can be expressed in terms of the Fox-Wright hypergeometric function of a single variable. Special cases of this integral lead to similar representations given in (\ref{e51}) and (\ref{e52}) for an integral considered by Ramanujan. Two infinite sums involving the Fox-Wright function are evaluated in closed form.

\vspace{0.6cm}

\noindent{\bf Acknowledgement} One of the authors (S.A.D.) was funded by the University Grants Commission of India through the award of a Dr. D.S. Kothari Post-doctoral Fellowship (Grant number F.4-2/2006(BSR)/MA/20-21/0061).

\end{document}